\documentclass{amsart}
\usepackage{amssymb}

\newtheorem{theorem}{Theorem}[section]
\newtheorem{lemma}[theorem]{Lemma}
\newtheorem{proposition}{Proposition}[section]   
\newtheorem{conjecture}{Conjecture}[section]

\theoremstyle{definition}
\newtheorem{definition}[theorem]{Definition}

\theoremstyle{remark}
\newtheorem{remark}[theorem]{Remark}

\numberwithin{equation}{section}

\renewcommand \H{U_q(\tilde{\mathfrak h})}
\renewcommand\a{\alpha}

\newcommand\ga{\gamma}
\newcommand{\Xa}[2]{x^{#1}_{i}(#2)}

\newcommand\ep{\epsilon}

\begin{document}

%Topmatter 

%One author
\title[New twisted quantum current algebras]
    {New twisted quantum current algebras
    %\\<title line 2>
    }
\author{Naihuan Jing}
\address{Department of Mathematics\\
   North Carolina State University\\ 
   Raleigh, NC 27695-8205\\
   U.S.A.} 
\email{jing@math.ncsu.edu}
\thanks{Research supported in part by NSA grant
MDA904-97-1-0062 and NSF grant DMS-9701755 at MSRI}

% Secondary: <subject>}
\date{January 15, 1999; revised on April 27, 1999}

%End topmatter

\begin{abstract}
We introduce a twisted quantum affine algebra associated
to each simply laced finite dimensional simple Lie algebra.
This new algebra is a Hopf algebra with a Drinfeld-type comultiplication. We obtain this algebra by considering its vertex representation. The vertex representation quantizes
the twisted vertex operators of Lepowsky-Wilson
and Frenkel-Lepowsky-Meurman. We also introduce a twisted quantum loop
algebra for the Kac-Moody case and give its level one representation.
  \end{abstract}

\maketitle

\section{Introduction}\label{S:1}

Let $\mathfrak g$ be the complex finite dimensional simple
Lie algebra of simply laced type. Let $\alpha_1, \cdots,
\alpha_l$ be the fixed set of simple roots. The associated
standard form is then defined by
\begin{equation}
(\a_i |\a_j)=a_{ij},
\end{equation}
where $A=(a_{ij})$ is the Cartan matrix.

Let $q$ be a generic complex number (i.e. not a root of unity). We define the $q$-integer $[n]$ by
\begin{equation*}
[n]=\frac{q^n-q^{-n}}{q-q^{-1}}.
\end{equation*}
Other $q$-numbers in this paper are defined accordingly using the same base. Thus the $q$-factorial $[n]!$ is defined by  $[n]!=\prod_{i=1}^n[i]$. 

In the theory of quantum affine algebras \cite{Dr, FJ}
the $q$-deformed Heisenberg algebra 
$U_q(\widehat{\mathfrak h})$ inside $U_q(\widehat{\mathfrak g})$ is an associative algebra
generated by $a_i(m)$ ($m\in \mathbb Z^{\times}$) and the
central element $\ga$ subject to the following relations
\begin{align}\label{E:1.1}
[a_i(m), a_j(n)]&=\delta_{m, -n} \frac{[(\a_i|\a_j)m]}{m} 
\frac{{\ga}^m-\ga^{-m}}{q-q^{-1}} ,\\
[a_i(m), \ga]&=0.
\end{align}

In \cite{FJ} we used this $q$-deformed Heisenberg algebra
to give the first explicit construction of level one irreducible modules of the quantum affine algebras
of simply laced types. Later the author \cite{J90}
generalized this construction to the twisted quantum affine algebras using some twisted $q$-Heisenberg algebras.
Recently in \cite{J98a} we introduced a quantum loop algebra for the
quantum Kac-Moody algebra by considering the bilinear
form as any symmetric form on an even lattice.

The purpose of this paper is to study an analogous Heisenberg
algebra and obtain a new twisted quantum affine algebra. 
Following the familiar twisting in conformal field theory \cite{FLM, L} we will restrict
the index set to be the set of odd integers in (\ref{E:1.1}).
With the help of this twisted Heisenberg algebra we will
construct a new class of twisted quantum algebras
and study their properties. Our construction also
follows the more general twisting discussed in \cite{J94}.

The classical twisted vertex operators played a special
role in the construction of the moonshine module $V^{\natural}$
\cite{FLM} for the Monster group. It also helped
build the twisted boson-fermion correspondence for the double covering group of $S_n$ \cite{J91a}. One can also
trace its appearance in various other
realizations of affine Lie algebras \cite{LP, DL}
and their quantum
analogs \cite{J96, DF}.

Part of our motivation is to construct a principal
bosonization of the quantum affine algebra $U_q(\widehat{sl}_2)$, which seems quite non-trivial
and is still unresolved.
The new class of twisted quantum algebras constructed in this paper is different from that of
the twisted quantum affine algebras given in \cite{J90}. It is
a common phenomenon that one classical object may have different
$q$-deformations. It is unclear what will be the relation between our new 
algebras and the 
untwisted quantum affine algebras, though at least in the case of $\widehat{sl}_2$ the author
conjectures that
 there exists a homomorphism from  $U_q(\widehat{sl}_2)$  into the
endomorphism ring of the Fock space constructed in section \ref{S:2}. We
suggest that the verification of this conjecture would provide a principal
construction for $U_q(\widehat{sl}_2)$. The difficulty seems to be the lack
of a Drinfled-Jimbo presentation for our twisted quantum current algebra. 

Replacing the underlying root lattice by the root lattice of a Kac-Moody
algebra, we also introduce a twisted quantum loop algebra associated to any
symmetrizable generalized Cartan matrix. The vertex operator representation of
this twisted Kac-Moody loop algebra uses some $q$-deformation of the rational
function $(\frac{1-x}{1+x})^r$. The commutation relations in this case involve
$q$-differences of the delta function. 

 The twisted current
algebras first appeared in our work \cite{J}  for
the simplest case of $A_1$. 
The algebra is similar to certain cases of Ding and Iohara's algebras
\cite{DI}, though we have relatively stricter relations. It appears that 
 a certain quotient of our algebra is also a quotient of their
algebra. One interesting phenomenon (see \cite{DI}) is that the Serre relations
involve with non-constant coefficients. In our case we replace the Serre
relations by commutation relations between two operators of same type (see
(\ref{E:3.11}-\ref{E:3.12}) ) as in the classical case \cite{FLM}.
 I thank Jintai Ding for helpful discussions and comments on the paper.

\section{Twisted vertex representations}\label{S:2}
Let $Q$ be the root lattice of the simple finite
dimensional Lie algebra $\mathfrak g$ of the simply laced type with the standard bilinear form $(\ \ |\ \ )$ such that
$(\a_i|\a_j)=a_{ij}$.

Let the $q$-deformed twisted Heisenberg algebra 
$U_q(\tilde{\mathfrak h})$ be the associative algebra
generated by $a_i(m)$ ($m\in 2\mathbb Z+1$) and the central element $\ga$ subject to the following relations
\begin{align}\label{E:2.1}
[a_i(m), a_j(n)]&=\delta_{m, -n} \frac{[(\a_i|\a_j)m]}{2m} 
\frac{{\ga}^m-\ga^{-m}}{q-q^{-1}} ,\\
[a_i(m), \ga]&=[a_i(m), \ga^{-1}]=0. \label{E:2.1a}
\end{align}

The algebra $\H$ has a canonical representation realized on 
$V=S(\tilde{\mathfrak h}^-)$, the space of 
symmetric polynomials in $a_{i}(-n), n\in 2\mathbb Z+1, n>0$.
This representation is actually the induced one from the trivial 
representation $1_{\mathbb C}$ with the central charge $\gamma=q$ of the subalgebra
$U_q(\mathfrak h^+)$ generated by $\gamma^{\pm 1}$ and 
$a_i(n), n\in 2\mathbb Z +1, n>0$. The explicit action
is given by
\begin{align}\notag
\gamma^{\pm 1} &= \text{the multiplication operator by $q^{\pm 1}$,}\\
a_i(-n)    &= \text{the multiplication operator by $a_i(-n)$,}\\
a_i(n)       &= \text{the annihilation operator by $a_i(n)$ subject to (\ref{E:2.1}),} \notag
\end{align}
for $n\in 2\mathbb Z+1 , n>0$.

We also denote by  $a_{i}(n) $ be the operator corresponding to $a_{i}(n)$, then they satisfy the relation (\ref{E:2.1}).

Let $\hat{Q}$ be the central extension of the
root lattice $Q$ such that
\begin{equation} \label{E:2.2}
1\longrightarrow \mathbb Z_2\longrightarrow \hat{Q}
\longrightarrow Q \longrightarrow 1
\end{equation}
with the commutator map 
\begin{equation}
aba^{-1}b^{-1}=(-1)^{(\a|\beta)},
\end{equation}
where we let $a$ be the preimage of $\a$. In the following we will denote the preimage of $\a_i$ by $a_i$. Let $T$ be 
any $\hat{Q}$-module such that as operators on $T$
\begin{equation}\label{E:2.3a}
a_ia_j=(-1)^{(\a_i|\a_j)}a_ja_i.
\end{equation}

Let $V_{Q}=S(\tilde{\mathfrak h}^-)\otimes T$.
We introduce the twisted vertex operators acting on
$V_{Q}$ by the following expressions:
\begin{align} \label{E:2.3}
E^{\pm}_{-}(\a_i,z) &= exp(\pm \sum^{\infty}_{n=1, odd} \frac{2 q^{\mp
n/2}}{[n]} a_i(-n) z^n), \\
E^{\pm}_{+}(\a_i, z) &= exp(\mp \sum^{\infty}_{n=1, odd} \frac{2 q^{\mp
n/2}}{[n]} a_i(n) z^{-n}) , \\
X_i^{\pm}(z)    &= E^{\pm}_{-}(\a_i, z)E^{\pm}_{+}(\a_i, z)a_i^{\pm 1}=\sum_{n\in\mathbb Z} X_i^{\pm}(n)z^{-n}.
\end{align}

We define the normal ordering as usual by moving the annihilation operators to the right of the creation operators. It is a routine calculation to obtain the 
operator product expansions for the twisted vertex operators.

\begin{proposition}\label{P:2.1}
\begin{align*}
X_i^{\pm}(z)X_j^{\pm}(w)&=:X_i^{\pm}(z)X_j^{\pm}(w):
\left\{\begin{array}{ll}
1  & \text{if $(\a_i|\a_j)=0$}\\
\frac{z-w}{z+w}\frac{z-q^{\pm 2}w}{z+q^{\pm 2}w} &
\text{if $i=j$}\\
\frac{z+q^{\mp 1}w}{z-q^{\mp 1}w} & \text{if $(\a_i|\a_j)=-1$}
\end{array}\right.     \\
X_i^{\pm}(z)X_j^{\mp}(w)&=:X_i^{\pm}(z)X_j^{\mp}(w):
\left\{\begin{array}{ll}
1& \text{if $(\a_i|\a_j)=0$}\\
\frac{z+qw}{z-qw}\frac{z+q^{-1}w}{z-q^{-1}w} &
\text{if $i=j$}\\
\frac{z-w}{z+w} & \text{if $(\a_i|\a_j)=-1$}
\end{array}\right.
\end{align*}

\end{proposition}
\medskip

In particular we have the following useful relations.

\begin{lemma} \label{P:2.2}
The operators 
$\Phi_i(zq^{-1/2})=:X_i^+(zq^{-1})X_i^-(z):$ and \newline
$\Psi_i(zq^{1/2})=:X_i^+(zq)X_i^-(z):$ are given by
\begin{align}
\Phi_i(z) &=exp((q^{-1}-q)\sum_{n>0, odd} a_i(-n)z^n)
=\sum_{n\geq 0} \Phi_i(-n)z^n ,\\
\Psi_i(z) &=exp((q-q^{-1})\sum_{n>0, odd} a_i(n)z^{-n})
=\sum_{n\geq 0}\Psi_i(n)z^{-n}. 
\end{align}
\end{lemma}
\medskip

These two exponential operators can generate a complete system of basis for
the $q$-Heisenberg algebra \cite{FJ, J90}. Let $G_{ij}(z)=\sum_{n\geq 0}G_nz^n$ be the Taylor series of 
the function 
\begin{equation} \label{E:g-fcn}
\frac{q^{(\a_i|\a_j)}z-1}{q^{(\a_i|\a_j)}z+1} \frac{z+q^{(\a_i|\a_j)}}{z-q^{(\a_i|\a_j)}}
\end{equation}
at $z=0$.

\begin{lemma} \label{P:2.3}
The relations (1.5) for the Heisenberg algebras
are equivalent to the following formal series one:
$$
\Phi_i(z)\Psi_j(w)=\Psi_j(w)\Phi_i(z)G_{ij}(\frac zw q^{-1})/G_{ij}(\frac zw q).
$$
\end{lemma}

\begin{proof} From (1.5) it follows that 
\begin{align*}
&\Phi_i(z)\Psi_j(w)\\ &=\Psi_j(w)\Phi_i(z)\,\exp\left\{\sum_{n\geq 0,odd}
\frac 2n(q^{(\a_i|\a_j)n}-q^{-(\a_i|\a_j)n})(q^n-q^{-n})\left(\frac zw\right)^n\right\}\\
&=\Psi_j(w)\Phi_i(z)G_{ij}(z/qw)/G_{ij}(qz/w),
\end{align*}
where we use the definition of $G_{ij}(z)$.
\end{proof}

It is easy to see that Lemma (\ref{P:2.3}) is equivalent
to the Heisenberg commutation relations (\ref{E:2.1}).
We give the complete relations satisfied by the twisted 
$q$-vertex operators as follows.

\begin{theorem} \label{P:2.4}
The twisted vertex operators obey 
the relations in Lemma \ref{P:2.3} and the following relations:
\begin{gather} \label{E:2.4}
[\Phi_i(z), \Phi_j(w)]=[\Psi_i(z), \Psi_j(w)]=0, \\
\Phi_i(z)X^{\pm}_j(w)\Phi_i(z)^{-1}=X^{\pm}_j(w)G_{ij}(q^{\mp 1/2}z/w)^{\pm 1}, \label{E:2.5}\\
\Psi_i(z)X^{\pm}_j(w)\Psi(z)^{-1}=X^{\pm}_j(w)G_{ij}(q^{\mp 1/2}w/z)^{\pm 1}, \label{E:2.6}\\
[X_i^{\pm}(z), X_j^{\pm}(w)]=[X_i^{+}(z), X_j^{-}(w)]=0, 
\qquad\text{$(\a_i|\a_j)=0$}   \label{E:2.7}\\
[X_i^+(z), X_j^-(w)]=2z:X_i^+(z)X_j^-(-z):\delta(-\frac wz), 
\qquad\text{$(\a_i|\a_j)=-1$} \label{E:2.8}\\
[X^+_i(z), X^-_i(w)]
=\frac {2(q+q^{-1})}{q-q^{-1}}
\left(\Psi_i(q^{-1/2}z)\delta(\frac wz q)-
\Phi_i(zq^{1/2})\delta(\frac wz q^{-1})\right),   \notag
\\
(z-q^{\mp 1}w)(z+q^{\pm 1}w)X^{\pm}_i(z)X_j^{\pm}(w)=
\label{E:2.9}\\
(q^{\mp 1}z-w)(q^{\pm 1}z+w)X^{\pm}_j(w)X_i^{\pm}(z), 
\qquad\text{$(\a_i|\a_j)=-1$}, \notag  \\
(z-q^{\pm 2}w)(z+q^{\mp 2}w)X^{\pm}_i(z)X^{\pm}_i(w)=
(q^{\pm 2}z-w)(q^{\mp 2}z+w) X^{\pm}_i(w) X^{\pm}_i(z) 
\label{E:2.10}
\\
%z^2(q-q^{-1})^2:X_i^{\pm}(z)^2:\delta(-\frac wz)
+2(q+q^{-1})^2\delta(-\frac wz) z^2,   \label{E:2.11}
\end{gather}
where $\delta(z)=\sum_{n\in \mathbb Z}z^n$ is the formal series of the 
$\delta$-function. 

\end{theorem}

\begin{proof} Relation (\ref{E:2.4}) is a trivial consequence of (\ref{E:2.1}).
Similar arguments as in the proof of Lemma (\ref{P:2.3}) show immediately
relations (\ref{E:2.5}-\ref{E:2.6}).

To prove the remaining relations we recall the basic property of the 
$\delta$-function: for a formal series $f(z)$ with $f(a)$ defined,
$$
f(z) \delta(\frac za)= f(a) \delta(\frac za).
$$

From Proposition \ref{P:2.1} it follows that for $(\a_i|\a_j)=2$
$$
%\multline
[X^+_i(z), X^-_i(w)] =
:X^+_i(z)X^-_i(w): \left( \frac{z+qw}{z-qw}\frac{z+q^{-1}w}{z-q^{-1}w}
-\frac{w+qz}{w-qw}\frac{w+q^{-1}z}{w-q^{-1}z}\right)
%\endmultline
$$
where the two rational functions in the parenthesis are formal series
in $w/z$ and  $z/w$ respectively.

Observe that 
$$
\frac{z+qw}{z-qw}\frac{z+q^{-1}w}{z-q^{-1}w}
=\frac{(z+qw)(z+q^{-1}w)}{zw} \sum_{n\geq 1}\frac {q^{n+1}-q^{-n-1}}
{q-q^{-1}}(\frac wz)^n
$$
as an identity of formal series. Therefore we have that 
\begin{align*} 
[X^+_i(z), X^-_i(w)] &=:X^+_i(z) X^-_i(w):
\frac{(z+qw)(z+q^{-1}w)}{zw(q-q^{-1})}
\left(\delta(\frac zw q^{-1}) -\delta(\frac zw q)\right) \\
&=\frac {2(q+q^{-1})}{q-q^{-1}}
\left(\Psi_i(q^{-1/2}z)\delta(\frac zw q^{-1})-
\Phi_i(zq^{1/2})\delta(\frac zw q)\right),
\end{align*}
where we used lemma \ref{P:2.2} and the property of the $\delta$-function. 
%The cases of (\ref{E:2.7}-\ref{E:2.8}) follow by a similar
%consideration.

Relation (\ref{E:2.9}) follows from similar consideration.
For $(\a_i|\a_j)=-1$ we have
\begin{align*}
&(z-q^{\mp 1}w)(z+q^{\pm 1}w)X^{\pm}_i(z)X_j^{\pm}(w)-
(q^{\mp 1}z-w)(q^{\pm 1}z+w)X^{\pm}_j(w)X_i^{\pm}(z)\\
&=:X^{\pm}_i(z)X_j^{\pm}(w):\left(
(z+q^{-1}w)-(w+q^{-1}z)\right)\\
&=0.
\end{align*}

Similarly we proceed to show (\ref{E:2.10}):
\begin{align*}
&(z-q^{\pm 2}w)(z+q^{\mp 2}w)X^{\pm}_i(z)X^{\pm}_i(w)-
(q^{\pm 2}z-w)(q^{\mp 2}z+w) X^{\pm}_i(w) X^{\pm}_i(z) \\
&=:X^{\pm}_i(z)X^{\pm}_i(w): (z-q^2w)(z-q^{-2}w)
\left(\frac{z-w}{z+w}+\frac{w-z}{w+z}\right)\\
&=:X^{\pm}_i(z)X^{\pm}_i(w):(z-q^2w)(z-q^{-2}w)(z-w)\delta(-\frac
zw)\\ 
&=2(1+q^2)(1+q^{-2})\delta(-\frac zw) z^2
%-z^2(q-q^{-1})^2:X^{\pm}_i(z)^2:,
\end{align*}
which completes the proof.
\end{proof}

\section{Twisted quantum algebras} \label{S:3}

We now introduce the new quantum affine algebra in this section. Let
$\mathfrak g$ be a complex simple Lie algebra with the standard bilinear
form $(\ \ |\ \ )$ on its root lattice $R$ as in section \ref{S:2}. 

\begin{definition}The algebra $U_q(\tilde{\mathfrak g})$ is an associative algebra generated by $a_{im}$, $x^{\pm}_{in}$, $\gamma=q^c$
($n\in \mathbb Z, m\in 2\mathbb Z+1$) with the following defining relations in terms of generating functions:
\begin{align}
\label{E:3.1}
x^{\pm}_i(z)&=\sum_{n\in\mathbb Z}x^{\pm}_{in}z^{-n}\\
\phi_i(z) &=exp((q^{-1}-q)\sum_{m>0, odd} a_{i,-m}z^m)
=\sum_{n\geq 0} \phi_{i,-n}z^n ,\label{E:3.2}\\
\psi_i(z) &=exp((q-q^{-1})\sum_{m>0, odd} a_{im}z^{-m})
=\sum_{n\geq 0}\psi_{in}z^{-n}. \label{E:3.3}
\end{align}
The defining relation are as follows.
\begin{gather} \label{E:3.4'}
[c, \phi_i(z)]=[c,\psi_i(z)]=[c, x_i^{\pm}(z)]=0,\\
\label{E:3.4}
[\phi_i(z), \phi_j(w)]=[\psi_i(z), \psi_j(w)]=0, \\
\phi_i(z)\psi_j(w)=\psi_j(w)\phi_i(z)G_{ij}(\frac zw q^{-c})/G_{ij}(\frac zw
q^c) ,  \label{E:3.5}\\ 
\phi_i(z)x^{\pm}_j(w)\phi_i(z)^{-1}=x^{\pm}_j(w)G_{ij}(q^{\mp c/2}z/w)^{\pm
1}, \label{E:3.6}\\
\psi_i(z)x^{\pm}_j(w)\psi(z)^{-1}=x^{\pm}_j(w)G_{ij}(q^{\mp c/2}w/z)^{\pm 1},
\label{E:3.7}\\ [x_i^{\pm}(z), x_j^{\pm}(w)]=[x_i^{+}(z), x_j^{-}(w)]=0, 
\qquad\text{$(\a_i|\a_j)=0$}   \label{E:3.8}\\ (z+w)[x_i^+(z), x_j^-(w)]=0, 
\qquad\text{$(\a_i|\a_j)=-1$} \label{E:3.9}\\ [x^+_i(z), x^-_i(w)]
=\frac {2(q+q^{-1})}{q-q^{-1}}
\left(\psi_i(q^{-c/2}z)\delta(\frac wz q^c)-
\phi_i(zq^{c/2})\delta(\frac wz q^{-c})\right),   \notag
\\
(z-q^{\mp 1}w)(z+q^{\pm 1}w)x^{\pm}_i(z)x_j^{\pm}(w)=
\label{E:3.10}\\
(q^{\mp 1}z-w)(q^{\pm 1}z+w)x^{\pm}_j(w)x_i^{\pm}(z), 
\qquad\text{$(\a_i|\a_j)=-1$}, \notag  \\
(z-q^{\pm 2}w)(z+q^{\mp 2}w)x^{\pm}_i(z)x^{\pm}_i(w)=
(q^{\pm 2}z-w)(q^{\mp 2}z+w) x^{\pm}_i(w) x^{\pm}_i(z) 
\label{E:3.11}\\
+2(q+q^{-1})^2\delta(-\frac wz)c. \label{E:3.12}
\end{gather}
where $G_{ij}$ is defined in (\ref{E:g-fcn}).
\end{definition}

 In the Drinfeld realization \cite{Dr} of the usual
quantum affine algebras there are several relations similar to our relations,
but with a different definition of $G_{ij}(z)$. Another obvious difference
is that we have more complicated commutation relations between operators
$x_i^{\pm}(z)$ of the same type.

 The
twisted algebra $U_q(\tilde{\mathfrak g})$ has a Hopf algebra structure. Let
$c_1=c\otimes 1$ and $c_2=1\otimes c$.

\begin{align}
\Delta(\Xa +z)&=\Xa +z\otimes 1+ \phi_{i}(q^{c_1/2}z)\otimes \Xa +{q^{c_1}z},\\
\Delta(\Xa -z)&=1\otimes \Xa +z + \Xa -{q^{c_2}z}\otimes \phi_{i}(q^{c_2/2}z),\\
\Delta(\phi_{i}(z))&=\phi_{i}(q^{-c_2/2}z)\otimes \phi_{i}(q^{c_1/2}z),\\
\Delta(\psi_{i}(z))&=\psi_{i}(q^{c_2/2}z)\otimes \psi_{i}(q^{-c_1/2}z),\\
\Delta(c)&=c\otimes 1+1\otimes c.
\end{align}

The antipode is given by:
\begin{align}
S(\Xa +z)&=-\phi_{i}({q^{-c/2}z})\Xa +{q^{-c}z},\\
S(\Xa -z)&=-\Xa -{q^{-c}z}\psi_{i}({q^{-c/2}z})^{-1},\\
S(\phi_{i}(z)) &=\phi_{i}(z)^{-1},\\
S(\psi_{i}(z)) &=\psi_{i}(z)^{-1},\\
S(c)&=-c.
\end{align}
The counit map is
\begin{align}
\ep(\Xa {\pm}z)&=0, \quad \ep(\phi_{i}(z))=1, \quad \ep(\psi_{i}(z))=1\\
\ep(q^c)&=1.
\end{align}

The relations (\ref{E:2.8}-\ref{E:2.10}) resolve the Serre relations for our
algebra. We remark that our algebra is different from the algebra of
Ding-Iohara \cite{DI} due to our relations (\ref{E:2.8}-\ref{E:2.10}).
If we replace the relation (\ref{E:2.10}) by the following weaker
one:
\begin{align*}
&(z+w)(z-q^{\pm 2}w)(z+q^{\mp 2}w)X^{\pm}_i(z)X^{\pm}_i(w)\\
&=(z+w)(q^{\pm 2}z-w)(q^{\mp 2}z+w) X^{\pm}_i(w) X^{\pm}_i(z). 
\end{align*}
Then this modified algebra has a common quotient with the
Ding-Iohara algebra in \cite{DI}.

\section{A twisted quantum Kac-Moody algebra and its level one representations} 
We can generalize the algebra to any even lattice. 
Let $Q$ be an arbitrary integral lattice generated by $\a_i$, $i=1, \cdots, l$,
equipped with a symmetric form $(\ \ |\ \ )$ such that
\begin{equation}\label{E:4.1}
(\a_i|\a_j)\in \mathbb Z, \qquad (\a_i|\a_j)=2.
\end{equation}
In the following we can view this lattice as that of a root lattice of
Kac-Moody algebra \cite{K1}. In this sense the following construction parallels
to the quantum loop algebra associated to the Kac-Moody algebra studied in
\cite{J98a}.

Let $U_q(\tilde{h}_{Q})$ be the associative algebra
generated by $a_i(m), m\in \mathbb Z$ and $\gamma^{\pm 1}$
with the defining relation (\ref{E:2.1}-\ref{E:2.1a}). Then we can consider the central extension of $Q$ by a finite group $<\kappa>$ similarly (cf. \cite{FLM}). If we use
$a_i$ to denote the preimage of $\a_i$ , then we have
\begin{equation}\label{E:4.2}
a_ia_j=(-1)^{(\a_i|\a_j)}a_ja_i.
\end{equation}

Let $T$ be any $\hat Q$-module such that
$\kappa=\omega$, a $o(\kappa)$th primitive root of unity.
The vertex space is defined accordingly:
\begin{equation}\label{E:4.3}
V_{Q}=S(\tilde{h}_{Q}^-)\otimes T.
\end{equation}

For each $\a_i$ we define the vertex operator by the same expression as in (\ref{E:2.3}). To compute the operator
product expansion we need the following $q$-analogs of  
series $(z-w)^r$ introduced in \cite{J96, J98a}.
For $r\in \mathbb C$ we define the $q$-analog to be in base $q^2$:
\begin{equation}
(1-z)_{q^{2}}^{r}=\frac{(q^{-r+1}z;q^{2})_{\infty}}
{(q^{r+1}z;q^{2})_{\infty}}=exp(-\sum_{n=1}^{\infty}\frac{[rn]}{n[n]}z^n),
\end{equation}
where $(b;q)_{\infty}=\prod_{n=0}^{\infty}(1-q^nb)$. The twisted $q$-analog is defined by
\begin{equation}
\left(\frac{1-z}{1+z}\right)_{q^{2}}^{r}=
\frac{(1-z)_{q^{2}}^{r}}{(1+z)_{q^{2}}^{r}}
=exp(-\sum_{n\in \mathbb Z^+_1}^{\infty}\frac{2[rn]}{n[n]}z^n),
\end{equation}

The following OPE can be calculated similarly as in section
\ref{S:2}.
\begin{proposition}\label{P:4.1}
\begin{align*}
X_i^{\pm}(z)X_j^{\pm}(w)&=:X_i^{\pm}(z)X_j^{\pm}(w):
\left(\frac{1-q^{\mp 1}w/z}{1+q^{\mp 1}w/z}\right)^{(\a_i|\a_j)}_{q^2}, \\
X_i^{\pm}(z)X_j^{\mp}(w)&=:X_i^{\pm}(z)X_j^{\mp}(w):
\left(\frac{1+w/z}{1-w/z}\right)^{(\a_i|\a_j)}_{q^2} .
\end{align*}
\end{proposition}
\medskip

Other commutations on the vertex representation space $V_{Q}$ are given as
follows. 

\begin{theorem} \label{P:4.4}
The twisted vertex operators satisfy the following relations:
\begin{gather} \label{E:4.4}
[\Phi_i(z), \Psi_j(w)]=[\Phi_i(z), \Psi_j(w)]=0, \\
\Phi_i(z)X^{\pm}_j(w)\Phi_i(z)^{-1}=X^{\pm}_j(w)G_{ij}(q^{\mp 1/2}z/w)^{\pm 1}, \label{E:4.5}\\
\Psi_i(z)X^{\pm}_j(w)\Psi(z)^{-1}=X^{\pm}_j(w)G_{ij}(q^{\mp 1/2}w/z)^{\pm 1}, \label{E:4.6}\\
[X_i^{\pm}(z), X_j^{\pm}(w)]=[X_i^{+}(z), X_j^{-}(w)]=0, 
\qquad\text{$(\a_i|\a_j)=0$,}   \label{E:4.7}\\
(z+w)^{-(\a_i|\a_j)}_{q^2}[X_i^+(z), X_j^-(w)]=0, 
\qquad\text{$(\a_i|\a_j)<0$} \label{E:4.8}\\
[X^+_i(z), X^-_i(w)]
=\frac {2(q+q^{-1})}{q-q^{-1}}
\left(\Psi_i(q^{-1/2}z)\delta(\frac wz q)-
\Phi_i(zq^{1/2})\delta(\frac wz q^{-1})\right),   \notag
\\
(z-w)_{q^2}^{-(\a_i|\a_j)-1}(z-q^{(\a_i|\a_j)}w)(z+q^{-(\a_i|\a_j)}w)X^{\pm}_i(z)X_j^{\pm}(w)=
\label{E:4.9}\\
(z-w)_{q^2}^{-(\a_i|\a_j)-1}(q^{(\a_i|\a_j)}z-w)(q^{-(\a_i|\a_i)}z+w)X^{\pm}_j(w)X_i^{\pm}(z), 
\qquad\text{$(\a_i|\a_j)<0$}, \notag  \\
(z-q^{\pm 2}w)(z+q^{\mp 2}w)X^{\pm}_i(z)X^{\pm}_i(w)=
(q^{\pm 2}z-w)(q^{\mp 2}z+w) X^{\pm}_i(w) X^{\pm}_i(z) 
\label{E:4.10}\\
+2(q+q^{-1})^2\delta(-\frac wz).  \label{E:4.11}
\end{gather}
\end{theorem}
\medskip

\begin{remark} Note that we use the factor $(z+w)^{-(\a_i|\a_j)}_{q^2}$ and
$(z-w)_{q^2}^{-(\a_i|\a_j)-1}$ to suppress poles in (\ref{E:4.8}) and
(\ref{E:4.10}) respectively. We  can cancel these two polynomials and rewrite
the commutation relations in terms of delta functions in a representation.  It
is always possible to refine the commutation relations as in the classical case
\cite{K2} and quantum case \cite{J98a}. In the quantum case one uses the
$q$-difference of the delta function $\delta(z-w)$ to expand zero expressions.
As in \cite{J98a} we have for $n\geq 0$ \begin{equation}
i_{z,w}(z-w)_{q^2}^{-n-1}-i_{w,z}(z-w)_{q^2}^{-n-1}
=\partial_{q,w}^{(n)}\delta(z-w),
\end{equation}
where $i_{z,w}$ means the expansion in the range
$|z|>>|w|$ and the $\partial_{q, w}^{(n)}$ denote the $q$-divided powers of
$n$th $q$-difference $\partial_q$:  \begin{equation}
\partial_q f(z)=\frac{f(qz)-f(q^{-1}z)}{(q-q^{-1})z}.
\end{equation}
 The delta function $\delta(z-w)$ means $z\delta(z/w)$. Then we can rewrite for
example the relation (\ref{E:4.9})  
\begin{align*}
&(z-q^{(\a_i|\a_j)}w)(z+q^{-(\a_i|\a_j)}w)X^{\pm}_i(z)X_j^{\pm}(w)\\
&\qquad -(q^{(\a_i|\a_j)}z-w)(q^{-(\a_i|\a_i)}z+w)X^{\pm}_j(w)X_i^{\pm}(z)\\
&=:X^{\pm}_i(z)X_j^{\pm}(w):(z+w)_{q^2}^{-(\a_i|\a_j)-1}
\partial_{q,w}^{(-(\a_i|\a_j))}\delta(z-w),
\end{align*}
where $(\a_i|\a_j)<0$. 
\end{remark}

\begin{definition} Let $\mathfrak g$ be a Kac-Moody algebra with the root
lattice $Q=\mathbb Z\alpha_1\oplus\cdots\oplus\mathbb Z\alpha_l$ and the
symmetric bilinear form $(\ \ | \ \ )$. The algebra $U_q(\widehat{\mathfrak
g})$ is an associative algebra with generators $a_{im}$, $x^{\pm}_{in}$, $\gamma=q^c$
($n\in \mathbb Z, m\in 2\mathbb Z+1$) and following defining relations:
\begin{gather} 
[c, \Phi_i(z)]=[c, \Psi_i(w)]=[c, X^{\pm}_i(z)]=0\\
[\Phi_i(z), \Psi_j(w)]=[\Phi_i(z), \Psi_j(w)]=0, \\
\Phi_i(z)X^{\pm}_j(w)\Phi_i(z)^{-1}=X^{\pm}_j(w)G_{ij}(q^{\mp c/2}z/w)^{\pm
1}, 
\\ \Psi_i(z)X^{\pm}_j(w)\Psi(z)^{-1}=X^{\pm}_j(w)G_{ij}(q^{\mp
c/2}w/z)^{\pm 1}, 
\\ [X_i^{\pm}(z), X_j^{\pm}(w)]=[X_i^{+}(z), X_j^{-}(w)]=0, 
\qquad\text{$(\a_i|\a_j)=0$,}   \\
(z+w)^{-(\a_i|\a_j)}_{q^2}[X_i^+(z), X_j^-(w)]=0, 
\qquad\text{$(\a_i|\a_j)<0$} \\
[X^+_i(z), X^-_i(w)]
=\frac {2(q+q^{-1})}{q-q^{-1}}
\left(\Psi_i(q^{-c/2}z)\delta(\frac wz q^c)-
\Phi_i(zq^{c/2})\delta(\frac wz q^{-c})\right),   \notag
\\
(z-w)_{q^2}^{-(\a_i|\a_j)-1}(z-q^{(\a_i|\a_j)}w)(z+q^{-(\a_i|\a_j)}w)X^{\pm}_i(z)X_j^{\pm}(w)=
\\
(z-w)_{q^2}^{-(\a_i|\a_j)-1}(q^{(\a_i|\a_j)}z-w)(q^{-(\a_i|\a_i)}z+w)X^{\pm}_j(w)X_i^{\pm}(z), 
\qquad\text{$(\a_i|\a_j)<0$}, \notag  \\
(z-q^{\pm 2}w)(z+q^{\mp 2}w)X^{\pm}_i(z)X^{\pm}_i(w)=
(q^{\pm 2}z-w)(q^{\mp 2}z+w) X^{\pm}_i(w) X^{\pm}_i(z) 
\\
+2(q+q^{-1})^2\delta(-\frac wz)c.  
\end{gather}
where the generating functions are defined by (\ref{E:3.2}--\ref{E:3.3}).
\end{definition}

\medskip

We close this section with the following statement.

\begin{conjecture} 
The quantum affine algebra $U_q(\hat{sl_2})$
is isomorphic to  the twisted algebra
 generated by $a_{1m} (m\in 2\mathbb Z+1),
\gamma, x_{1n}^{\pm} (n\in \mathbb Z)$ subject to the relations
(\ref{E:3.4'}-\ref{E:3.12}) associated to $A_1$. 
\end{conjecture}
At $c=1$ this statement is true by specializing $q\to 1$
\cite{LW},
however we do not know the exact isomorphism at the 
quantum level. The exact isomorphism will solve the so-called principal vertex construction 
of $U_q(\widehat{sl}_2)$.

\end{document}